\documentclass[10pt]{article}
\usepackage{amssymb}
\usepackage{amsmath}
\usepackage{amsfonts}
\usepackage{graphicx}

\usepackage{latexsym}
\title{Copulas in three dimensions with prescribed correlations}
\author{Luc Devroye and G\'erard Letac}

\begin{document}

\maketitle

\begin{abstract} Given an arbitrary three-dimensional correlation matrix,
we prove that there exists a three-dimensional joint distribution for the random variable $(X,Y,Z)$ such that $X$,$Y$ and $Z$ are identically distributed with beta
distribution $\beta_{k,k}(dx)$ on $(0,1)$ if $k\geq 1/2$. 
This implies that any correlation structure can be attained for
three-dimensional copulas.
\end{abstract}


\section{Introduction} 

The simulation community is quite interested in the computer generation
of identically distributed random vectors $(X_1, \ldots , X_n)$ with
prescribed marginal distribution ($\mu$) and correlation matrix ($R$).
We call a fixed distribution $\mu$ on the Borel sets of $\mathbb{R}$
\textit{$n$-universal} if for every possible correlation matrix,
there exists a joint distribution for $(X_1, \ldots , X_n)$
that achieves it, while respecting the marginal distribution condition.
We call $\mu$ \textit{universal}  if it is $n$-universal for all $n$.
For example, the standard normal law is universal: just decompose the
(positive semi-definite) correlation matrix $R$ into its Choleski form
$S \times S^t$. Then verify that if $X$ is a column vector of $n$
i.i.d.\ normal random variables, then $SX$ has covariance matrix
$\mathbb{E} (SXX^tS^t ) = S \times S^t = R$.
If $Y$ is any random variable not identically zero, and $N$ is standard
normal, then the normal scale mixture law of $YN$ is universal---just
check that $YX$ has the same correlation matrix as $X$, where $X$ is as above.
It would be of general interest to characterize all universal
distributions.

The blossoming field of copulas (see, e.g., Nelsen, 2006) is
largely concerned with similar issues, but until now it was mainly
interested in $n=2$ (see, e.g., Genest and MacKay (1986)
for some early work), and, by convention, in the uniform marginal law $\mu$.
Marginal distributions without atoms can be mapped to a uniform law by the probability
integral transform, and back with the inverse probability integral transform,
hence the central role of the uniform law. These transforms generally alter
the correlation matrix, but some transforms, such as between the normal
and the uniform do not alter correlations a lot (see, e.g., Falk, 1999).
This means that the universality problem has to be tackled for each marginal
$\mu$ separately.

For example, because of its universality, it is convenient to
start with a random vector $(X_1, \ldots , X_n)$ with normal marginals
as constructed above, achieving a certain correlation matrix $R$.
Then $(\Phi (X_1), \ldots , \Phi (X_n))$ is a random vector with
uniform $[0,1]$ marginals. Falk (1999) discusses the merits of this approach
by noting that the maximal deviation between correlation coefficients
before and after is at most 0.0181. 
One could attempt to start with a different correlation matrix $R'$ for
the normal random  vector in the hope of obtaining the right correlation
matrix $R$ for $(\Phi (X_1), \ldots , \Phi (X_n))$.
However, for $n \ge 3$, this strategy is doomed to fail for some $R$,
no matter how hard one tries in the construction of $R'$: this
will be shown in Section 4.

Approximative solutions abound in the literature---some of these
are surveyed in Devroye (1986). New approximations are being
developed regularly, see, e.g., Headrick (2009). 
Discrete laws where the marginals are all Bernoulli have received
particular attention (Emrich and Piedmonte (1991), Lee (1993)), but
for any fixed marginal structure, the possible values of the $2^n$ ``free
joint probabilities'' form a polytope, and thus, the region of allowable
correlation matrix coefficients forms a polytope as well. As we will see below, 
for $n \ge 3$, the region of all allowable correlation matrix coefficients
is convex but is  not a polytope (see figure below for $n=3$). 
In fact, any marginal law $\mu$ with a finite number of atoms is
not $n$-universal for $n \ge 3$ for this reason.

It is well-known that the uniform law and many other laws are 2-universal.
The question begs whether the uniform law is universal. We do
not have the answer to this, but the purpose of this note is to
lift the veil carefully, and to reveal that the uniform law
is 3-universal.  In fact, we will show that all symmetric beta laws
of parameter $k \ge 1/2$ are 3-universal:
\begin{equation}\label{BB}
\beta_{k,k}(dx)=\frac{\Gamma(2k)}{\Gamma^2(k)}x^{k-1}(1-x)^{k-1}\textbf{1}_{(0,1)}(x)dx.
\end{equation}
Arguing as we did above for the normal law, this implies that
all symmetric unimodal densities are 3-universal---just observe that
by Khinchine's theorem, each symmetric unimodal random variable can be written as $Y U$,
where $U$ is uniform $[-1,1]$ and $Y$ is arbitrary and independent of $U$.

For $p,q,r\in \mathbb{R}$ consider the $(3,3)$ symmetric matrix
\begin{equation}\label{RC}R=\left[\begin{array}{ccc}1&r&q\\r&1&p\\q&p&1\end{array}\right]\end{equation}We denote by $\mathcal{R}_n$ the set of semipositive definite matrices  of order $n$ with 
unit diagonal elements.
By computing the principal minors one sees that $R$ defined by (\ref{RC}) is in $ \mathcal{R}_3$ if and only if the numbers $1-p^2,1-q^2,1-r^2$ and 
\begin{equation}\label{DG} 
\Delta=\Delta(p,q,r)=\det R=1-p^2-q^2-r^2+2pqr\end{equation}
are nonnegative. 

Recall that if $X$ and $Y$ are uniformly distributed on $(0,1)$ an abundant literature calls the joint distribution of $(X,Y)$ a copula.
For this reason,
let us call the distribution of $(X_1,\ldots,X_n)$
an $n$-dimensional copula when $X_1,\ldots,X_n$ are uniformly distributed on $(0,1)$.


\section{Facts about $\mathcal{R}_n$ and $n$-dimensional copulas}

Trivially the set $\mathcal{R}_n$ is a closed convex subset of the linear space of symmetric matrices of order $n$. However the set of its extreme points is not easy to grasp and its characterization given  by Ycart (1986) is difficult to handle for $n\geq 4$.
It makes a sharp distinction with the cone $\mathcal{P}_n$ of positive 
semi-definite symmetric matrices of order $n$,
whose extremal lines are generated by matrices of rank one as an easy consequence of the spectral theorem for symmetric matrices. In the case of $\mathcal{R}_n$ there  are only $2^{n-1}$ matrices of rank one. These matrices  are
$(\epsilon_i\epsilon _j)_{1\leq i,j\leq n}$ where $\epsilon_j=\pm 1$ (they are extreme points of $\mathcal{R}_n$ and they are the correlation matrices of $X(\epsilon_1,\ldots,\epsilon_n)$ when $X$ is a one-dimensional random variable with a second moment). Another notable point is  the fact that positive definite correlation matrices are inner points of  $\mathcal{R}_n$ and cannot be extreme points: there are no extreme points of rank $n$.
Finally,
for our purposes,
the most important result of Ycart (1986) is his characterization of the extreme points of $\mathcal{R}_n$ which are matrices of rank 1 and 2. They are the matrices of the form  
$$(\cos (\alpha_i-\alpha_j))_{1\leq i,j\leq n}.
$$
where $\alpha_1,\ldots,\alpha_n$ are arbitrary numbers. 

\begin{figure}[h!]
  \centering
      \includegraphics[width=0.5\textwidth]{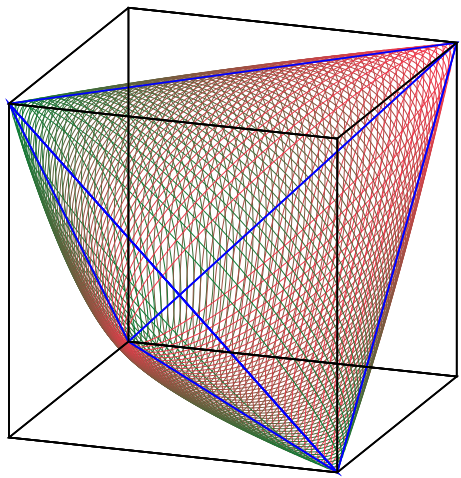}
  \caption{The space of the three off-diagonal correlation coefficients of a correlation matrix is a convex subset of $[0,1]^3$. Strictly contained in it is the simplex formed by the four
extremal points $(1,1,1), (1,-1,-1), (-1,-1,1), (-1,1,-1)$. As a first step,
we will give a construction for the three-dimensional [uniform] copula that
achieves correlation equal to the point $(-1/2, -1/2, -1/2)$ on the bulge of the surface. This surface is invariant by a 4-element group of rotations 
generated by the three symmetries with respect to the three axes.}
\end{figure}

The consideration of the extreme points of $\mathcal{R}_n$ is justified by the following observation: suppose that $(X_1,\ldots,X_n)\sim \mu$ and $(X'_1,\ldots,X'_n)\sim \mu'$,
suppose that the real random variables $X_1,\ldots,X_n,X'_1,\ldots,X'_n$ are identically distributed and have second moment and denote by $R(\mu)$ and $R(\mu')$ the respective correlation matrices of $\mu$ and $\mu'$.
Let $\lambda\in [0,1]$ and $(X''_1,\ldots,X''_n)\sim \lambda\mu+(1-\lambda)\mu'$.
Then obviously $X''_j\sim X_1$ and the correlation matrix of the mixing 
satisfies \begin{equation}\label{RCC}R(\lambda\mu+(1-\lambda)\mu')=\lambda R(\mu)+(1-\lambda)R(\mu').\end{equation}
Therefore in order to prove that for any $R\in \mathcal{R}_n$ there exists a $n$-dimensional copula $\mu$ such that $R=R(\mu)$ enough is to prove it for all the cases where $R$ is an extreme point of $\mathcal{R}_n$.
For from the Caratheodory theorem,
since $\mathcal{R}_n$ is immersed in an affine space of dimension $d=n(n-1)/2$ then for any $R\in \mathcal{R}_n$ there exist $d+1$ extreme points $R_0,\ldots,R_d$ of $\mathcal{R}_n$ and $d+1$ nonnegative numbers $\lambda_0,\ldots \lambda_d$ of sum 1 such that $R=\sum_{j=0}^d\lambda_jR_j$.
If we have been able to find for each extreme point $R_j$ an $n$-dimensional copula $\mu_j$ such that $R_j=R(\mu_j)$ then from (\ref{RCC}) we get
$$
R=R \left( \sum_{j=0}^d\lambda_j\mu_j \right) .
$$
Needless to say these remarks extend to the case where the uniform distribution on $(0,1)$ is replaced by a distribution with finite second moments.

Let us comment on the cases $n=2$ and $n=3$.

\vspace{4mm}\noindent \textsc{The case $n=2$.} Clearly the elements of $\mathcal{R}_2$ have the form $R(r)=\left[\begin{array}{cc}1&r\\r&1\end{array}\right]$ where $-1\leq r\leq 1$  and  the two extreme points are $R(1)$ and $R(-1)$ since
$$R(r)=\frac{1+r}{2}R(1)+\frac{1-r}{2}R(-1)
$$
As explained above this leads to an immediate solution to the problem of finding a two-dimensional copula with correlation matrix $R(r)$.

If $X$ is uniformly distributed,
denote by $\mu_1$ the distribution of $(X,X)$,
by $\mu_{-1}$ the distribution of $(X,-X)$,
and by $\mu_r$ the mixing
$\mu_r=\frac{1+r}{2}\mu_1+\frac{1-r}{2}\mu_{-1}$ leading to $R(r)=R(\mu_r)$.
Later in Section 4 we shall give another useful way of
using Gaussian variables to design a two-dimensional copula with a given correlation matrix.

\vspace{4mm}\noindent \textsc{The case $n=3$.} As seen before,
the extreme points of $\mathcal{R}_3$ have rank 1 or 2 and  they have the form (\ref{RC}) with $p=\cos(\alpha_2-\alpha_3),\ q=\cos(\alpha_3-\alpha_1),\ r=\cos(\alpha_1-\alpha_2)$.
For convenience we rather write $a=\alpha_2-\alpha_3,
b=\alpha_3-\alpha_1,
c=\alpha_1-\alpha_2$,
with 
$$a+b+c\equiv 0\ \mathrm{mod}\ 2\pi.
$$
With this constraint we denote
\begin{equation}\label{RCT}
R(a,b,c)=\left[\begin{array}{ccc}1&\cos c&\cos b\\\cos c&1&\cos a\\\cos b&\cos a&1\end{array}\right] .
\end{equation}
For instance, we get the four matrices of rank one  with the choices 
$$(a,b,c)=(0,0,0),
\ (0,\pi,\pi),
\  (\pi,0,\pi),
\ (\pi,\pi,0).
$$
It is a good exercise to check 
$\Delta(\cos a,\cos b,\cos c)=0$ with the notation (\ref{DG}). Note that $R(a,b,c)$ has rank 2 if and only if $(\sin a,\sin b,\sin c)\neq (0,0,0)$.


\section{Construction}

\vspace{4mm}\noindent \textsc{ A geometric construction in a particular case.} The aim of this section is to build a three-dimensional copula with the extremal correlation matrix $R(a,b,c)$ defined by (\ref{RCT}) with $a+b+c\equiv 0\ \mathrm{mod}\ 2\pi$.
As explained in Section 2,
this enables us to find a three-dimensional copula with an arbitrary correlation matrix $R$.
We shall be able to do this even  by replacing the uniform distribution $\beta_{1,1}$ by the beta distribution $\beta_{k,k}$ defined by (\ref{BB}) with $k\geq 1/2$.
Since the solution is relatively complicated,
it is desirable to work first for a special choice of parameters which has been for us a path toward the general case. This particular case is $p=q=r=-1/2$ in the notation (\ref{RC}) or $a=b=c=2\pi/3$ in the notation (\ref{RCT}). Since the correlation matrix is the covariance matrix of some affine transformations of the initial random variables,
we rather consider three random variables $(X,Y,Z)$ which are uniformly distributed on $(-\sqrt{3},\sqrt{3})$---thus centered with variance 1.
We want to construct the distribution of $(X,Y,Z)$ such that 
$$\mathbb{E}(XY)=\mathbb{E}(YZ)=\mathbb{E}(ZX)=-1/2,\ \mathbb{E}(X^2)= \mathbb{E}(Y^2)= \mathbb{E}(Z^2)=1.
$$
This implies that  $\mathbb{E}((X+Y+Z)^2)=0$ and $X+Y+Z=0$ almost surely. The intersection of the plane $X+Y+Z=0$ and the cube $\{|X|,|Y|,|Z|\leq \sqrt{3}\}$ is a hexagon. The largest disk $D$ contained in this hexagon is centered at 0 and has radius $3/\sqrt{2}$ since it is tangent to the sides of the hexagon at  six points
$$\pm(\sqrt{3}/2,
\sqrt{3}/2,-\sqrt{3}),
\ \pm(\sqrt{3}/2,-\sqrt{3},
\sqrt{3}/2),
\ \pm(-\sqrt{3},\sqrt{3}/2,
\sqrt{3}/2).
\ 
$$
We now define an appropriate (and unique) distribution for  $(X,Y,Z)$ on $D$ such that it is invariant by the rotations of the disk and such that the projections $(X,Y,Z)\mapsto X$ are uniform on $(-\sqrt{3},\sqrt{3})$.

An orthonormal basis of the plane $X+Y+Z=0$ is the pair of vectors $f_1=(0,1/\sqrt{2},-1/\sqrt{2})$ and $f_2=(\sqrt{2}/\sqrt{3},-1/\sqrt{6},-1/\sqrt{6})$.
Thus a distribution on $D$ which is invariant by rotation is the distribution of 
$$\frac{3}{\sqrt{2}}R(\cos \Theta f_1+\sin \Theta f_2) 
$$
where $\Theta$ is uniform on $[0,2\pi)$ and is independent of the random variable $R\in (0,1)$. This leads to
$$X=\sqrt{3}R\sin \Theta,\ Y=-\sqrt{3}R\sin (\Theta-\frac{\pi}{3}),\ Z=-\sqrt{3}R\sin (\Theta+\frac{\pi}{3}).
$$ 
Because $\Theta$ is uniform,
$X,$ $Y$ and $Z$ are identically distributed. Now we take $R$ with distribution $\frac{r}{ \sqrt{1-r^2}}\textbf{1}_{0,1)}(r)dr$ and we show that the distribution of $R|\sin \Theta|$ is uniform on $(0,1)$ by computing its Mellin transform. For $s>0$ we have
\begin{eqnarray}\label{FF}\mathbb{E}(R^s)\mathbb{E}(|\sin \Theta|^s)&=&\int_0^1r^s\frac{r}{ \sqrt{1-r^2}}dr \times \frac{2}{\pi}\int_{0}^{\pi/2}(\sin \theta)^{2\frac{s+1}{2}-1}d\theta\\&=&\frac{1}{2}B(1+\frac{s}{2},\frac{1}{2})\times\frac{1}{\pi}B(\frac{1}{2}+\frac{s}{2},\frac{1}{2}) =\frac{1}{1+s} . \nonumber \  \square\end{eqnarray}
\medskip

\vspace{4mm}\noindent \textsc{Comments.} We have just given  an analytic proof of a theorem due to Archimedes,
which says that if you project the uniform distribution on the three-dimensional sphere $S_2$ onto a diameter,
you get the uniform distribution on the diameter. Consider the three-dimensional sphere $S$ constructed from the disk $D$ above (We mean: having center zero and radius equal to the radius of $D$ namely $3/\sqrt{2})$.
Put the uniform distribution on $S$,
project it orthogonally on $D$: this projection is actually the distribution of $(X,Y,Z)$ above. 

\begin{figure}[h!]
  \centering
      \includegraphics[width=\textwidth]{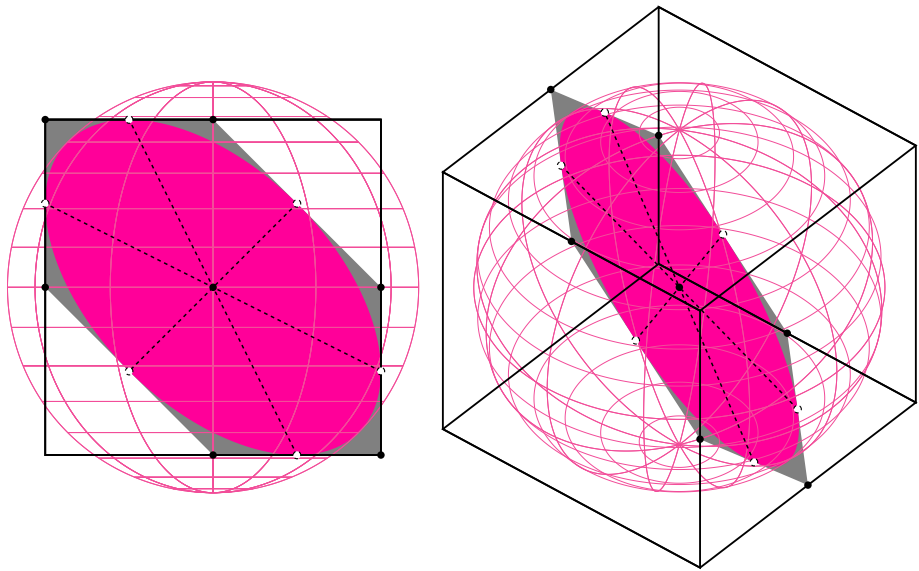}
  \caption{Illustration of our construction. First take a point uniformly on the surface of the ball. Project it to the plane shown (so that it falls in the circle). The three coordinates of that point
are each uniformly distributed on $[0,1]$. Also, the correlation structure is given by
$\mathbb{E}(XY)=\mathbb{E}(YZ)=\mathbb{E}(ZX)=-1/2,\ \mathbb{E}(X^2)= \mathbb{E}(Y^2)= \mathbb{E}(Z^2)=1.$}
\end{figure}

\vspace{4mm}\noindent \textsc{The general case.}
Without loss of generality we assume that $R(a,b,c)$ defined by (\ref{RCT}) has rank two,
that is, $(\sin a,\sin b,\sin c)\neq (0,0,0)$.
In the sequel we assume without loss of generality that $\sin c\neq 0$.

Let us observe that  $R(a,b,c)$ is the correlation matrix of the centered random variable $(X,Y,Z)$ if and only if 
\begin{equation}\label{PP}X\sin a+Y\sin b+Z\sin c=0.\end{equation}
To see this, observe that 
\begin{eqnarray*}
&&\mathbb{E} ((X\sin a+Y\sin b+Z\sin c)^2) \\
&&=\sin^2 a+\sin^2 b+\sin^2 c+2(\cos c\sin a\sin b+\cos b\sin a\sin c+\cos a\sin b\sin c)\\
&&=\sin^2 (b+c)+\sin^2 b+\sin^2 c\\
&&\qquad  +2(-\cos c\sin (b+c)\sin b-\cos b\sin b+c)\sin c+\cos (b+c)\sin b\sin c)\\
&&=0.
\end{eqnarray*}
We imitate the previous particular case as follows: we watch the intersection of the plane $P$ defined by (\ref{PP}) with the cube $(-1,1)^3$ (which is slightly more convenient for the general case than the cube$(-\sqrt{3},\sqrt{3})^3$ that we have used before). This intersection is a hexagon---and not a lozenge---for the following reason:

\vspace{4mm}\noindent \textbf{Lemma 3.1.} The plane $\alpha X+\beta Y+\gamma Z=0$ defines a plane $P$ such that the intersection with $(-1,1)^3$ is a hexagon if and only if there exists a triangle with sides $|\alpha|,|\beta|,|\gamma|,$ namely
$$
|\alpha|<|\beta|+|\gamma|,\ |\beta|< |\alpha|+|\gamma|,\ |\gamma|<|\alpha|+|\beta|.
$$ 

\vspace{4mm}\noindent \textbf{Proof.} $\Rightarrow:$ Assume that $|\alpha|\geq |\beta|+|\gamma|$.
Then the intersections of $P$ with the four lines $Y=\epsilon=\pm 1,\ Z=\eta=\pm 1$ are given by the four points
$$
\left( -\epsilon\frac{\beta}{\alpha}-\eta\frac{\gamma}{\alpha},\epsilon,\eta \right) ,
$$
which are the vertices of a lozenge. Since $-1<-\epsilon\frac{\beta}{\alpha}-\eta\frac{\gamma}{\alpha}<1$ this lozenge in contained in the cube. The converse is similar. $\square$
\medskip

With $a+b+c\equiv 0 \ \mathrm{mod}\,
2\pi$ we now apply the lemma to $\alpha=\sin a,\ \beta=\sin b,\ \gamma=\sin c$.
Observe that
$$|\sin a|=|\sin (b+c)|=|\sin b\cos c+\sin c \cos b|\leq |\sin b|+|\sin c|
$$
and the two other inequalities are similar: therefore the intersection of $P$ with the cube is a hexagon $H$.

We now construct the unique ellipse $E$ inscribed in $H$.
For this we introduce for $z\in (-1,1)$ the ellipse $E_z$ in the $(x,y)$ plane defined by
$$E_z=\{(x,y)\ ; \Delta(x,y,z)=0\}
$$
where $\Delta$ is defined by (\ref{DG}). This ellipse is inscribed in the square $[-1,1]^2$ and the four points of tangency with the square are $\pm (z,1)$ and $\pm(1,z)$.
We denote by $U_z$ the closed convex hull of $E_z$. 

\vspace{4mm}\noindent \textbf{Lemma 3.2.} The projection of the ellipse $E$  by $(x,y,z)\mapsto (x,y)$ is $E_{\cos c}$

\vspace{4mm}\noindent \textbf{Proof.} Denote by $E'$ the intersection of the plane $P$ with the cylinder 
$$\{(x,y,z)\ ; \Delta(x,y,\cos c)=0\}.
$$
We want to prove that $E'$ is inscribed in $H$ and thus that $E=E'$. Since the four points of tangency of $E_{\cos c}$ with the square $[-1,1]^2$ are 
$\pm (\cos c,1)$ and $\pm(1,\cos c),$ observe that the  unique points of the four lines $\pm (\cos c,1,z)$ and $\pm(1,\cos c,z)$ which are located on $P$ are exactly $A_{\pm}=\pm(1,\cos c,\cos b),$ and  $B_{\pm}=\pm(\cos c,1,\cos a)$.
By construction these four points are points of tangency of $E'$ with $H$.
Now observe that the two points $C_{\pm}=\pm(\cos b,\cos a,1)$ belong to $E'$: This comes from the fact that $\Delta(\cos b,\cos a,\cos c)=0$ (since the matrix $R(a,b,c)$ of (\ref{RCC}) is singular) and thus $C_{\pm}$ are on the cylinder; the second fact is that 
$\sin a \cos b+\sin b\cos a+\sin c=0$ and this implies that $C_{\pm}$ are on $P$.
Trivially also $C_{\pm}$ are in $H$. 
For showing that $C_{\pm}$ are points of tangency of $E'$ to $H$ we have to prove that $(x,y,z)\in E'$ implies that $|z|\leq 1$.
Here is a not too elegant proof: we have to deduce $|z|\leq 1$ from the two equalities
$ \Delta(x,y,\cos c)=0$  and $x\sin a +y\sin b+z\sin c=0$ or
$$x^2+y^2-2xy\cos c=\sin ^2 c,\ -z=\frac{\sin a}{\sin c}x+\frac{\sin b}{\sin c}y.
$$
We make the change of variable 
$x=\frac{1}{\sqrt{2}}(X+Y),\ y=\frac{1}{\sqrt{2}}(X-Y)$, which leads to
$$
\frac{X^2}{1-\cos c}+\frac{Y^2}{1+\cos c}=1 .
$$
We therefore write $X=\sqrt{1-\cos c}\cos \theta,\ Y=\sqrt{1+\cos c}\sin \theta$ and we get
$$
x=\sin(\theta+\frac{c}{2}),\ y=\sin(\theta-\frac{c}{2}),\ z=-\cos \tau\cos \theta-\sin \tau \sin \theta=-\cos (\theta-\tau) ,
$$
where 
\begin{equation}\label{TE}
\cos \tau=\left(\frac{\sin a}{\sin c}+\frac{\sin b}{\sin c}\right)\cos\frac{c}{2},\ \sin \tau=\left(\frac{\sin a}{\sin c}-\frac{\sin b}{\sin c}\right)\sin\frac{c}{2}
\end{equation}
(we skip the proof that the squares of these two  quantities add to one and that $\tau$ does exist). Finally $|z|=|\cos (\tau -\theta)|\leq 1$. 
To conclude the proof of the lemma,
we observe that $E'$ has six distinct points with $H$,
$A_{\pm},\  B_{\pm},\ C_{\pm}$.  
This implies that $E'$ is the unique ellipse inscribed in $H$ and therefore $E=E'$.
This ends the proof of the lemma. $\square$ 
\medskip

Because of the symmetry between $(a,b,c)$ and the previous lemma the two other projections of $E$ by $(x,y,z)\mapsto (x,z)$ and $(x,y,z)\mapsto (y,z)$ have respectively the equations $\Delta(x,\cos b,z)=0$ and $\Delta(\cos a,y,z)=0$.

We now construct a probability $\mu_k$ on the convex hull $U$ of $E$ such the three margins of $\mu_k$ are the probabilities 
$$\nu_k(dx)=\frac{2^{1-2k}\Gamma(2k)}{\Gamma^2(k)}(1-x^2)^{k-1}\textbf{1}_{(-1,1)}(x)dx
$$
where $k\geq 1/2$.
Note that $\nu_k$ is nothing but the image of $\beta_k$ by the map $x\mapsto 2x-1$.
For this construction of $\mu_k$ we use the following parameterization of $U$ by $r\in (0,1]$ and $\theta\in [0,2\pi)$ inherited from the calculations of the previous lemma,
where $\tau$ is defined by (\ref{TE})
\begin{eqnarray}
x&=&r\sin (\theta+\frac{c}{2})\nonumber\\
y&=&r\sin( \theta-\frac{c}{2})\label{PRE}\\
z&=&-r\cos (\theta-\tau) .     \nonumber
\end{eqnarray}
On the space $(0,1]\times [0,2\pi)$ we define the probability for $k>1/2$ 
$$\rho_k(dr,d\theta)=(2k-1)(1-r^2)^{k-\frac{3}{2}}rdr\times \frac{d\theta}{2\pi}
$$
as well has 
$$
\rho_{1/2}(dr,d\theta)=\delta_1(dr)\times \frac{d\theta}{2\pi}
$$
and on the space $U$ we define the probability $\mu_k$ as the image of $\rho_k$ by the map $(r,\theta)\mapsto (x,y,z)$
defined by the equalities (\ref{PRE}). If $(X,Y,Z)\sim \mu_k$ we claim now that $X\sim Y\sim Z\sim \nu_k$.
The proof  by Mellin transform is quite analogous to (\ref{FF}). 
For $k>1/2$ if $\Theta$  is uniform on $[0,2\pi)$  and is independent of $R\sim (2k-1)(1-r^2)^{k-\frac{3}{2}}r\textbf{1}_{(0,1)}(r)dr$ then $ (R,\Theta)\sim \rho_k$.
If $R=1$ the same is true for $k=1/2$.

We show that the distribution of $R|\sin \Theta|$  is the distribution of $|T|$ when $T\sim \nu_k$ by computing its Mellin transform. 
For $s>0$ we have
\begin{eqnarray*}\label{FG}\mathbb{E}(R^s)\mathbb{E}(|\sin \Theta|^s)&=&\frac{2k-1}{2}\int_0^1(1-r^2)^{k-\frac{3}{2}}r^{s+1}dr \times \frac{2}{\pi}\int_{0}^{\pi/2}(\sin \theta)^{2\frac{s+1}{2}-1}d\theta\\&=&(k-\frac{1}{2})B(1+\frac{s}{2},k-\frac{1}{2})\times\frac{1}{\pi}B(\frac{1}{2}+\frac{s}{2},\frac{1}{2})\nonumber\\&=& \frac{1}{\sqrt{\pi}}\frac{\Gamma(k+\frac{1}{2})\Gamma(\frac{s+1}{2})}{\Gamma(k+\frac{s+1}{2})}=\frac{2^{-2k-1}}{B(k,k)}B(\frac{s+1}{2},k)\nonumber\\&=&\frac{2^{-2k}}{B(k,k)}\int_0^1(1-t^2)^{k-1}t^sdt=\mathbb{E}(|T|^s) .
\end{eqnarray*}
Since $R\sin \Theta$ and $T$ are symmetric random variables we get $R\sin \Theta\sim T$.
Finally
$$X=R\cos (\Theta-\frac{c}{2})\sim Y=R\cos (\Theta-\frac{c}{2})\sim Z=-R\cos (\Theta-\tau)\sim R\sin \Theta\sim T.
$$
The Mellin transforms of $|T|$ and $R$  have shown that
$$
\mathbb{E}(X^2)=\mathbb{E}(Y^2)=\mathbb{E}(Z^2)=\frac{1}{2k+1},\ \mathbb{E}(R^2)=\frac{2}{2k+1}.
$$
Therefore,
$$
\mathbb{E}(XY)=\frac{\mathbb{E}(R^2)}{2}\mathbb{E}(2\sin (\Theta-\frac{c}{2})\sin (\Theta+\frac{c}{2}))=\frac{\cos c}{2k+1}.
$$
This shows that the correlation between $X$ and $Y$ is $\cos c$.
to get the correlation between $X$ and $Z$ we use 
$Z=-\frac{\sin a}{\sin c}X-\frac{\sin b}{\sin c}Y$ and the computation already done of $\mathbb{E}(XY)$ and of  $\mathbb{E}(X^2)$. We get easily that the correlation of $(X,Z)$ is $ \cos b$.
Similarly  the correlation between $Y$ and $Z$
is $\cos a$,
and this achieves the proof of 
$R(\mu_k)=R(a,b,c)$.

\vspace{4mm}\noindent \textsc{Comments about the two-dimensional marginals.} 
In the general case we do no longer have for $\mu_k$ the generalization of the beautiful interpretation of 
$\mu_1$ in terms of the Archimedes theorem.
However if $(X,Y,Z)\sim \mu_k$ the joint distributions of $(X,Y)$,
$(Y,Z)$ and $(Z,X)$ already appear in the literature. 
The distribution $\phi_{k,\cos c}(dx,dy)$ of $(X,Y) $  is concentrated on the convex hull $U_{\cos c}$ of the ellipse $E_{\cos c}$ when $k>1/2$ and  is concentrated on the ellipse $E_{\cos c}$ for $k=1/2$.
For $k>1/2$ we have since $\Delta(x,y,\cos c)=(1-r^2)\sin^2 c$
\begin{equation*}\label {PHI}
\phi_{k,\cos c}(dx,dy)=
\frac{2k-1}{2\pi} |\sin c|^{\frac{1}{2}-k}\Delta(x,y,\cos c)^{k-\frac{3}{2}}\textbf{1}_{U_{\cos c}}(x,y)\,
dxdy .
\end{equation*}
This distribution $\phi_{k,\cos c}$ appears as a Lancaster distribution for the pair $(\nu_k,\nu_k)$ 
More specifically consider the sequence $(Q_n)_{n=0}^{\infty}$ of the orthonormal polynomials for the weight $\nu_k.$
 Thus $Q_n$ is the Jacobi polynomial $P_n^{k-1,k-1}$ 
normalized such that $\int_{-1}^1Q^2_n(x)\nu_k(dx)=1.$ For $1/2<k$
denote $$K(x,y,z)=\sum_{n=0}^{\infty}\frac{Q_n(x)Q_n(y)Q_n(z)}{Q_n(1)}.$$ This series converges if $|x|,|y|,|z|<1$
and its sum is zero when $(x,y)$ is not in the interior $U_z$ of the ellipse $E_z.$ With this notation we have 
$$\phi_{k,\cos c}(dx,dy)=K(x,y,\cos c)\nu_k(dx)\nu_k(dy).$$ This result is essentially due to Gasper (1971). See  Koudou (1995) and (1996) and Letac (2009) for details. 
Needless to say,
the distributions of $(Y,Z)$ and $(Z,X)$
are $\phi_{k,\cos a}$ and $\phi_{k,\cos b}$.


\section{The use of Gaussian variables for building a $n$-dimensional copula  with given correlation.}

We start from the simplest idea: if 
$$\Phi(x)=\frac{1}{\sqrt{2\pi}}\int_{-\infty}^xe^{-\frac{u^2}{2}}du
$$
and if $X\sim N(0,1)$  we have $\Phi (X)$ uniform on $(0,1)$. Let us denote
\begin{equation}\label{TTT}T(x)=T(x)=2\sqrt{3}(\Phi(x)-1/2).
\end{equation}
Then $T(X)$ is uniform on $(-\sqrt{3},\sqrt{3})$ with mean 0 and variance 1. Assume now that $(X_1,\ldots,X_n)$ is Gaussian $N(0,R)$ with covariance $R\in \mathcal{R}_n$.
Denote by $R^*$ the covariance matrix of $(T(X_1),\ldots,T(X_n))$ and observe that $R^*$
is also the correlation matrix of $(T(X_1),\ldots,T(X_n))$ and of the $n$-dimensional copula $(\Phi(X_1),\ldots,\Phi(X_n))$. [Copulas constructed in this manner are sometimes
called Gaussian copulas.]
In this section we compute $R^*$ as a function of $R$ and we examine the image $\mathcal{R}^*_n$ of $\mathcal{R}_n$ into itself of the function $R\mapsto R^*$. As we are going to see $\mathcal{R}^*_n$ is strictly smaller than $\mathcal{R}_n$ for $n\geq 3$.
To compute $R^*$, we need the following result.

\vspace{4mm}\noindent \vspace{4mm}\noindent \textbf{Proposition 4.1.} 
Let $(X,Y)$ be a centered Gaussian variable of $\mathbb{R}^2$ with covariance matrix $\left[\begin{array}{cc}1&r\\r&1\end{array}\right]$.
Then 
\begin{eqnarray}\label{PHI}
\mathbb{E}(\Phi(X)\Phi(Y))&=&\frac{1}{2}-\frac{1}{2\pi}\arg \cos\frac{r}{2} , \\
\label{THI}
\mathbb{E}(T(X)T(Y))&=&3-\frac{6}{\pi}\arg \cos\frac{r}{2} .
\end{eqnarray}
The proof of (\ref{PHI}) can be done by brute force and the computation of a four-dimensional integral. We rather going to obtain Proposition 4.1 in a more interesting way after the following result.

\vspace{4mm}\noindent \textbf{Theorem 4.2.} Let $(X,Y)$ be a centered Gaussian variable of $\mathbb{R}^2$ with covariance matrix $\left[\begin{array}{cc}1&r\\r&1\end{array}\right]$ and let $f$ be a real measurable function such that $\mathbb{E}_r(f(X))=0$ and $\mathbb{E}_r(f(X)^2)=1$.
Consider the Hermite polynomials $(H_n)_{n=0}^{\infty}$ defined by the generating function
$$e^{xt-\frac{t^2}{2}}=\sum_{n=0}^{\infty}H_n(x)\frac{t^n}{n!}
$$
and the expansion in orthogonal functions 
$$f(x)=\sum_{n=1}^{\infty}a_n\frac{H_n(x)}{\sqrt{n!}}.
$$
Then for all $-1\leq r\leq 1$ we have $\sum_{n=1}^{\infty}a_n^2=1$ and 
\begin{equation}\label{AAA}\mathbb{E}(f(X)f(Y))=\sum_{n=1}^{\infty}a_n^2r^n.\end{equation}

\vspace{4mm}\noindent \textbf{Proof.} 
Let us compute 
$$
\mathbb{E} \left( e^{Xt-\frac{t^2}{2}}e^{Ys-\frac{s^2}{2}} \right)
= \sum_{n=0}^{\infty}\sum_{m=0}^{\infty}\frac{t^n}{n!}\frac{s^m}{m!}\mathbb{E}(H_n(X)H_m(Y)).
$$
For this,
write $r=\cos \alpha$ with $0\leq \alpha\leq \pi$.
If $X,Z$ are independent centered real Gaussian random variables with variance 1,
then  $Y=X\cos \alpha +Z\sin \alpha $ is centered with variance 1,
$(X,Y)$ is Gaussian and $\mathbb{E}(XY)=\cos \alpha$. 
Therefore a simple calculation  gives 
$$
\mathbb{E} \left( e^{Xt-\frac{t^2}{2}}e^{Ys-\frac{s^2}{2}} \right)
=e^{ts\cos \alpha}.
$$
This shows that 
$\mathbb{E}(H_n(X)H_m(Y))=0$ if $n\neq m$ and that $\mathbb{E}(H_n(X)H_n(Y))=n!\cos^n \alpha$.
From this we get the result. $\square$
\medskip

\vspace{4mm}\noindent \textbf{Corollary 4.3.} Let $p_n\geq 0$ such that $\sum_{n=1}^{\infty}p_n=1$ and consider the generating function $g(r)=\sum_{n=1}^{\infty}p_n r^n$.
Let $R=(r_{ij})_{1\leq i,j\leq d}$ be in $\mathcal{R}_n$.
Then $R^*=(g(r_{ij}))_{1\leq i,j\leq d}$ is the covariance and correlation matrix of the random variable 
$(f(X_1),\ldots,f(X_d))$ where $(X_1,\ldots,X_d)$ is centered Gaussian with covariance $R$ and where 
$$f(x)=\sum_{n=1}^{\infty}\epsilon _n\sqrt{p_n}\frac{H_n(x)}{\sqrt{n!}}
$$
with  $\epsilon_n=\pm 1$.

\vspace{4mm}\noindent \textbf{Comment.} Note that many functions $f$ can give the same covariance $R^*$ for $(f(X_1),\ldots,f(X_d)),$ by taking arbitrary signs in the sequence $(\epsilon_n)_{n\geq 1}$ above.  It is worthwhile mentioning that Theorem  4.2 is easily extended to $p$ variables in the following sense. Consider 
$$
f(x^{(1)},\ldots,x^{(p)})=\sum_{n_1,\ldots,n_p}a_{n_1,\ldots,n_p}\frac{H_{n_1}(x^{(1)})\ldots H_{n_p}(x^{(p)})}{\sqrt{n_1!\ldots n_p!}},
$$
and assume that $a_{0,\ldots,0}=0$ and that $\sum_{n_1,\ldots,n_p}a^2_{n_1,\ldots,n_p}=1$.
Define
$$
g(r^{(1)},\ldots,r^{(p)})=\sum_{n_1,\ldots,n_p}a^2_{n_1,\ldots,n_p}(r^{(1)})^{n_1}\ldots (r^{(p)})^{n_p}.
$$
Let $X^{(k)}=(X_1^{(k)},\ldots,X_n^{(k)})\sim N(0,R^{(k)})$
for $k=1,\ldots,p$ be independent Gaussian variables of $\mathbb{R}^n$ such that the covariance $R^{(k)}=(r^{(k)}_{ij})_{1\leq i,j\leq n}$ is a correlation matrix.
Define $Y_j=f(X_j^{(1)},\ldots, X_j^{(p)})$.
Then the covariance matrix of $(Y_1,\ldots, Y_n) $ is the correlation matrix 
$$
(g(r^{(1)}_{ij},\ldots,r^{(p)}_{ij}))_{1\leq i,j\leq n}.
$$
\medskip

\vspace{4mm}\noindent \textbf{Proof of Proposition 4.1.} 
We apply Theorem 4.2 to the function $f=T$ defined by (\ref{TTT}). 
For this we have to compute 
$$
\frac{a_n}{\sqrt{n!}}=\mathbb{E} \left( T(X)\frac{H_n(X)}{n!} \right).
$$
Note that this is zero for even $n$  since $H_n$ and $T$ respectively even and odd functions. Thus we have to compute $p_{2n+1}> 0$ and $\epsilon_{2n+1}=\pm 1$ such that 
$$
\epsilon_{2n+1}\frac{\sqrt{p_{2n+1}}}{\sqrt{(2n+1)!}}=\mathbb{E} \left( T(X)\frac{H_{2n+1}(X)}{(2n+1)!} \right).
$$
To this purpose we watch the coefficient of $t^n$ in the power expansion of 
$$
\mathbb{E} \left( T(X)e^{Xt-\frac{t^2}{2}} \right).
$$
For this we need 
$$
\mathbb{E} \left( \Phi(X)e^{Xt-\frac{t^2}{2}} \right)
=\Phi \left( \frac{t}{\sqrt{2}} \right)
=\frac{1}{2}+\frac{1}{2\sqrt{\pi}}\sum_{n=0}^{\infty}\frac{(-1)^n}{4^nn!}\frac{t^{2n+1}}{2n+1},
$$
and
$$
\mathbb{E} \left( T(X)e^{Xt-\frac{t^2}{2}} \right)
=\sqrt{\frac{3}{\pi}}\sum_{n=0}^{\infty}\frac{(-1)^n}{4^nn!}\frac{t^{2n+1}}{2n+1}.
$$
Therefore ,
$$
\epsilon_{2n+1}\frac{\sqrt{p_{2n+1}}}{\sqrt{(2n+1)!}}=\sqrt{\frac{3}{\pi}}\frac{(-1)^n}{4^nn!}\frac{1}{2n+1},
$$
which shows that $\epsilon_{2n+1}=(-1)^n$. 
To finish the proof we apply (\ref{AAA}) to $a_{2n+1}=\sqrt{p_{2n+1}}$ and we get
$$
\mathbb{E}(T(X)T(Y)=\frac{3}{\pi}\sum_{n=0}^{\infty}(\frac{1}{2})_n\frac{1}{4^nn!}\frac{r^{2n+1}}{2n+1}=3-\frac{6}{\pi}\arg \cos\frac{r}{2},
$$
the last equality being easily checked. Of course (\ref{PHI}) is deduced from (\ref{THI})$.\square$
\medskip

We now prove that $\mathcal{R}^*_3$ is strictly smaller than $\mathcal{R}_3$.
For this we observe that
$$
r^*=3-\frac{6}{\pi} \arccos \frac{r}{2}\Leftrightarrow r=2\sin \frac{\pi r^*}{6}.
$$
Consider the matrix
$$
R_p=\left[\begin{array}{ccc}1&p&p\\p&1&p\\p&p&1\end{array}\right].
$$
Since $\det R_p=(1-p)^2(1+2p)$,
 the matrix $R_p$ is in $\mathcal{R}^*_3$ if and only if $-1/2\leq p\leq 1$.
However $2\sin \frac{\pi}{12}<-1/2$ therefore $R_{2\sin \frac{\pi}{12}}$ cannot be a correlation matrix.
This shows that there is no Gaussian variable $(X,Y,Z)$ such that the correlation matrix of 
$(\Phi(X),\Phi(Y),\Phi(Z))$ is $R_{-1/2}$.
To see that $\mathcal{R}^*_n$ is strictly smaller than $\mathcal{R}_n$ for $n\geq 4$,
observe that that the block matrix $\mathrm{diag}(R_{-1/2},I_{n-3})$ is in $\mathcal{R}_n$ and not in $\mathcal{R}^*_n$.
 
\section{Acknowledgments}
The authors thank Hakan Demirtas for helpful discussions. G\'erard Letac thanks  Sapienza Universit\`a di Roma for its generous support during the preparation of this paper. 

\section{References}

\vspace{4mm}\noindent 
\textsc{Devroye,
L.} (1986)
\textit{``Non-Uniform Random Variate Generation.''}
Springer-Verlag,
New York.

\vspace{4mm}\noindent 
\textsc{Emrich,
M. J.,
and Piedmonte,
M.R.} (1991)
``A method for generating high-dimensional multivariate binary variates'' 
 \textit{Amer. Statist.},
\textbf{45},
302-304.

\vspace{4mm}\noindent 
\textsc{Falk,
M.} (1999)
``A simple approach to the generation of uniformly distributed random variables with prescribed correlations''
\textit{Comm. Statist. Simulation Comp.},
\textbf{28},
785-791.

\vspace{4mm}\noindent 
\textsc{Gasper, G.} (1971)
``Banach algebra for Jacobi series and positivity of a kernel''
\textit{Ann. of Math.},
\textbf{95}, 261-280.

\vspace{4mm}\noindent 
\textsc{Genest, C.\ and MacKay, J.} (1986)
``The joy of copulas: bivariate distributions with uniform marginals''
 \textit{Amer. Statist.},
\textbf{40},
280-283.

\vspace{4mm}\noindent
\textsc{Headrick, T. C.} (2009)
\textit{``Statistical Simulation: Power Method Polynomials and other Transformations.''}
Chapman \& Hall / CRC Press, Boca Raton, FL.

\vspace{4mm}\noindent 
\textsc{Koudou,
A. E.} (1995),
\textit{``Probl\`emes de marges et familles exponentielles naturelles.''}
Th\`ese,
Universit\'e Paul Sabatier,
Toulouse.
 
\vspace{4mm}\noindent 
\textsc{Koudou,
A. E.} (1996) ``Probabilit\'es de Lancaster'' \textit{Expositiones Math.} \textbf{14}, 247-275.

\vspace{4mm}\noindent 
\textsc{Lee,
A. J.} (1993)
``Generating random binary deviates having fixed marginal distributions and specified degrees of association''
 \textit{Amer. Statist}.,
\textbf{47},
209-215.

\vspace{4mm}\noindent 
\textsc{Letac,
G.} (2008) 
``Lancaster probabilities and Gibbs sampling'' 
{\it  Statistical  Science},
\textbf{23},
187-191.

\vspace{4mm}\noindent 
\textsc{Nelsen,
R. B.} (2006)
\textit{``An Introduction to Copulas.''}
Springer-Verlag,
Berlin.

\vspace{4mm}\noindent 
\textsc{Ycart,
B.} (1985) 
``Extreme points in convex sets of symmetric matrices''
\textit{Proceedings of the American Mathematical Society}
\textbf{95},
issue 4,
607-612. 

\end{document}